\documentclass[12pt]{amsart}
\usepackage{amsfonts}

\theoremstyle{definition}

\theoremstyle{remark}

\numberwithin{equation}{section}

\begin{document}
\title{On multidimensional F. Riesz's "Rising Sun"
Lemma}

\author{A.A. Korenovskyy}
\address{Department of Mathematical Analysis, IMEM,
National University of Odessa, Dvoryanskaya, 2, 65026 Odessa,
Ukraine}
\email{anakor@paco.net}
\thanks{The work of the first author is partially supported by
the France-Ukraine program of the scientific collaboration
"DNIPRO"}

\author{A.K. Lerner}
\address{Department of Mathematics,
Bar-Ilan University, 52900 Ramat Gan, Israel}
\email{aklerner@netvision.net.il}

\author{A.M. Stokolos}
\address{Department of Mathematics,
University of Connecticut, U-9, Storrs, CT 06268, USA}
\email{stokolos@math.uconn.edu}

\subjclass{Primary 42B25}

\date{}
\keywords{"Rising sun" lemma, dyadic property, differential basis}
\commby{}

\begin{abstract}
A multidimensional version of the Riesz rising sun lemma is proved
by means of a generalized dyadic process.
\end{abstract}

\maketitle

The "rising sun" lemma of F. Riesz \cite{Riesz} is very important
in the differentiation theory, in the theory of the
one-dimensional Hardy-Littlewood maximal function (see
\cite{GMS,Riesz1}), and, as was mentioned by E.M. Stein\cite{Ste},
"...had implicitly played a key role in the earlier treatment of
the Hilbert transform".

\vskip 1mm \noindent {\bf Lemma (F. Riesz (cf. \cite{Klemes})).}
{\it Let $f$ be a summable function on some interval $I_0\subset
{\mathbb R}$, and suppose $|I_0|^{-1}\int_{I_0}f(x)\,dx\le A$.
Then there is a finite or countable set of pairwise disjoint
subintervals $I_j\subset I_0$ such that
$\left|I_j\right|^{-1}\int_{I_j} f(x)\,dx=A \ (j=1,2,\dots)$, and
$f(x)\le A$ for almost all $x\in I_0\setminus\left(\bigcup_{j\ge
1}I_j \right)$. } \vskip 1mm

All known proofs of this lemma (and its variants) are based on the
topological structure of the real line. In particular the lemma
does not rely on Lebesgue measure and works for arbitrary
absolutely continuous measure. This fact can be used, for example,
to show that the one-dimensional Hardy-Littlewood maximal function
with respect to any measure is of weak type $(1,1)$. As is well
known, for $n\ge 2$ the situation is quite different \cite{Sjo}.

A classical multidimensional substitute of Riesz's lemma is the
Calder\'on-Zygmund lemma \cite{CZ}, which probably one of the most
important propositions of harmonic analysis.

\vskip 1mm \noindent {\bf Lemma (Calder\'on-Zygmund).} {\it Let
$f\ge 0$ be a summable function on some cube $Q_0\subset {\mathbb
R}^n$, and suppose $|Q_0|^{-1}\int_{Q_0}f(x)\,dx\le A$. Then there
is a finite or countable set of pairwise disjoint subcubes
$Q_j\subset Q_0$ such that $A<\left|Q_j\right|^{-1}\int_{Q_j}
f(x)\,dx\le 2^nA$ $(j=1,2,\dots),$ and $f(x)\le A$ for almost all
$x\in Q_0\setminus\left(\bigcup_{j\ge 1}Q_j\right)$. } \vskip 1mm

In the one-dimensional case an "interval" and a "cube" are the
same, thus Riesz's lemma is a more precise result than the
Calder\'on-Zygmund lemma. This makes the use of Riesz's lemma is
more preferable for obtaining some sharp results. For instance,
sharp estimates of non-increasing rearrangements of functions from
$BMO$~\cite{Klemes}, of functions satisfying Gehring, Muckenhoupt,
and Gurov-Reshetnyak type conditions \cite{Ko1,Ko2,Ko3} were
obtained by means of Riesz's lemma.

A natural question arises whether one can improve the
Calder\'on-Zygmund lemma for $n\ge 2$, in other words, whether one
can get a multidimensional variant of Riesz's lemma where by an "
interval" we mean a "cube". Unfortunately, simple examples show
that this is not the case. Indeed, it suffices to consider
$f(x)=1- \chi_{[1/2,1]^2}(x), Q_0=[0,1]^2$, and a number close
to~$1$, for instance, $A=7/8$.

It is well known that the geometry of rectangular intervals
$I=\prod_{i=1}^n\left[a_i,b_i\right)$ (which we will simply call
rectangles) is much richer in comparison with cubes. Therefore, it
is either very surprising or very natural that rectangles play a
main role in an actual multidimensional analogue of Riesz's lemma,
and this is the result which we shall prove.

\vskip 1mm \noindent {\bf Lemma.} {\it Let $I_0$ be a rectangle in
$ {\mathbb R}^n$, and let $\mu$ be an absolutely continuous Borel
measure on $I_0$. Let $f\in L_\mu(I_0)$ and
$\frac{1}{\mu(I_0)}\int_{I_0}f(x)d\mu\le A$. Then there is a
finite or countable set of pairwise disjoint rectangles
$I_{j}\subset I_{0}$ with
$(\mu\left(I_j\right))^{-1}\int_{I_j}f(x)d\mu=A\,$
$(j=1,2,\dots),$ and $f(x)\le A$ for $\mu$- almost all points
$x\in I_0\setminus \left( \bigcup _{j\ge 1}I_{j}\right)$.} \vskip
1mm

Before proving this lemma, several remarks are in order. \vskip
1mm \noindent {\it Remark 1.} Crucial steps in this direction were
made by the first author who suggested some weak version of
Riesz's lemma in the multidimensional case \cite{Ko4,Ko5}. \vskip
1mm

\vskip 1mm \noindent {\it Remark 2.} In the two-dimensional case
the lemma contains implicitly in Besicovitch's paper \cite[Lemma
1]{Bes}. However, an attempt to extend the method of that paper to
$n\ge 3$ allows us only to get the lemma with an additional
restriction on grows $f$, namely, $f$ must belong to
$L(\log^+L)^{n-2}$. \vskip 1mm

\noindent {\it Remark 3.} As we mentioned above, known proofs of
Riesz's lemma are essentially based on the order structure of the
real line. The proof of the Calder\'on-Zygmund lemma is different,
namely, it relies on conception of dyadic cubes and on
stopping-time argument. Nevertheless, both proofs are unified by
the fact that the required intervals or cubes are form the level
set $\{x\in P_0:\widetilde Mf(x)>A\}$. In case of Riesz's lemma,
$P_0=I_0$ and $\widetilde Mf$ is the one-sided (left or right)
Hardy-Littlewood maximal function. In case of the
Calder\'on-Zygmund lemma, $P_0=Q_0$ and $\widetilde Mf$ is the
dyadic maximal function. Our proof does not rely on the strong
maximal function. More precisely, our main idea is a generalized
dyadic process, that is, for any number $A$ we construct an
individual differential basis of rectangles which has nice
covering properties like the basis of dyadic cubes. \vskip 1mm

\noindent {\it Remark 4.} In the case $n=1$ we get a new simple
proof of Riesz's lemma itself.

\begin{proof}[Proof of the lemma.]
If $f_{I_0}=A$, then the lemma is obvious. Assume that
$f_{I_0}<A$. Divide $I_0$ into two rectangles by the
$(n-1)$-dimensional hyperplane passing through the middle of the
largest side. The mean value of $f$ over at least one of the
partial rectangles must be less then $A$. If the mean value of $f$
over both rectangles is less then $A$, we continue the division of
each of them. Otherwise we obtain that the mean value of $f$ over
one of them is less then $A$, while over the another one is bigger
then $A$. In this case we move the hyperplane parallel to itself
until we get two rectangles such that the mean value of $f$ over
one of them is exactly equal to $A$, while over the another one is
less then $A$ (here we take in mind the absolutely continuity of
$\mu$). We put the bigger rectangle (over which the mean value is
equal to $A$) into the family $\left\{I_j\right\}$, while the
smaller one will be further divided. Thus we have described the
division argument.

We will continue the above process applying each time the division
argument to the rectangle which should be divided. As a result we
will get two family of rectangles: a finite or countable family of
pairwise disjoint rectangles $\left\{I_j\right\}$, such that
$f_{I_j}=A\ (j\ge 1)$, and a countable family of rectangles
$\left\{J_j\right\}$ for which $f_{J_j}<~A\ (j\ge~1)$.

We now observe that the rectangles from $\left\{J_j\right\}$ have
the following "dyadic" property: if $J_i\cap J_k\not=\emptyset$,
then either $J_i\subset J_k$ or $J_k\subset J_i$ for all
$J_i,J_k\in \left\{J_j\right\}$. Since, by the division argument,
each time we divided the biggest side at least on half, we get
that for any $x\not\in \cup_j I_j$ there is a sequence of
rectangles $J_m\in \left\{J_j\right\}$ containing $x$ and such
that their diameters tend to zero. Therefore, rectangles
$\left\{J_j\right\}$ form a differential basis (see \cite[Ch. 2,\S
2]{Guz} for definitions) on the set $E\equiv
I_0\setminus\left(\bigcup_{j\ge 1}I_j\right)$ with nice covering
properties. In particular, clearly that "dyadic" property implies
the Vitaly covering property \cite[Ch. 1]{Guz}, and hence the
basis ${\mathcal B}=\left\{J_j\right\}$ differentiates $L^1(E)$.
By the standard differentiation argument we obtain $f(x)\le A$ for
$\mu$-almost all $x$ from $E$.
\end{proof}

\end{document}